\documentclass[11pt]{article}
\usepackage[cp1251]{inputenc}
\usepackage[T2A]{fontenc}
\usepackage[english]{babel}
\usepackage{blindtext}
\usepackage[margin=2cm,left=2cm,includefoot]{geometry}
\usepackage[mathcal]{eucal}
\usepackage{amsfonts}
\usepackage{amssymb}
\usepackage{ amssymb }
\usepackage{amsmath}
\usepackage{hyperref}
\usepackage{abstract}
\usepackage[symbol]{footmisc}

\usepackage{amsthm}
\usepackage{amscd}
\usepackage{mathrsfs}
\usepackage{ dsfont }
\usepackage{xcolor}
\usepackage{blindtext}

\usepackage{amsfonts}
\usepackage{amssymb}
\usepackage{ amssymb }
\usepackage{amsmath}
\usepackage{abstract}

\theoremstyle{plain}

\newtheorem{lem}{Lemma}
\newtheorem{foll}{Corollary}

\newtheorem{theor}{Theorem}

\newtheorem{prop}{Proposition}

\newenvironment{solve}{\begin{proof}[Proof]}{\end{proof}}

\newtheorem{thmx}{Theorem}

\begin{document}

\begin{center}
\textsc{\Large Derivative of the iterations of Minkowski question mark function in special points.\footnote[2]{This research is supported by the Russian Science Foundation under grant [19-11-00065] and performed in Khabarovsk Division of the Institute for Applied Mathematics, Far Eastern Branch, Russian Academy of Sciences.}}\\ 
\large Nikita Shulga
\end{center}

\begin{abstract}
\noindent
For the Minkowski question mark function $?(x)$ we consider derivative of the function $f_n(x) = \underbrace{?(?(...?}_\text{n times}(x)))$. Apart from obvious cases (rational numbers for example) it is non-trivial to find explicit examples of numbers $x$ for which  $f'_n(x)=0$. In this paper we present a set of irrational numbers, such that for every element $x_0$ of this set and for any $n\in\mathbb{Z}_+$ one has $f'_n(x_0)=0$.
\end{abstract}

\section{Introduction}
For $x\in [0,1]$ we consider its continued fraction expansion
$$x=[a_1,a_2,\ldots,a_n,\ldots]=  \cfrac{1}{a_1+\cfrac{1}{a_2+\cdots}}, \,\,\, a_j \in \mathbb{Z}_+$$
which is unique and infinite when $x\not\in\mathbb{Q}$  and finite for rational $x$. 
Each rational $x$ has just two representations
$$
x=[a_1,a_2,\ldots,a_{n-1}, a_n]  
\,\,\,\,\,\text{and}\,\,\,\,\,
x=[a_1,a_2,\ldots,a_{n-1}, a_n-1,1]
,
\,\,\,\, \text{where}\,\,\,\, a_n \ge 2
.
$$
\\
For irrational $x=[a_1,a_2,\ldots,a_n,\ldots]$   the formula
\begin{equation} \label{eq:1}
?(x)=\sum\limits_{k=1}^{\infty} \frac{(-1)^{k+1}}{2^{a_1+\ldots+a_k-1}}
\end{equation}
introduced by Denjoy \cite{Den,Den1} and Salem \cite{Salem}
may be considered as one of the equivalent definitions of the function $?(x)$, which was first introduced by Minkowski in 1904 in \cite{Mink}. If $x$ is rational, then the infinite series in (\ref{eq:1}) is replaced by a finite sum, so that 
\begin{equation}\label{rational}
?([a_1,\ldots,a_t+1])=?([a_1,\ldots,a_t,1])=\sum\limits_{k=1}^{t} \frac{(-1)^{k+1}}{2^{a_1+\ldots+a_k-1}}
\end{equation}
 and hence $?(x)$ is well-defined for rational numbers by formula \eqref{rational} too. Derivative of the Minkowski question mark function was studied by several authors, see for example \cite{MD}, \cite{dushistova}, \cite{Kan1}, \cite{Kan2}, \cite{paradis}. 

In this paper we will consider $n$-th iteration of the $?(x)$ function, i.e. the function
$$f_n(x):= \underbrace{?(?(...?}_\text{n times}(x))).$$
Iterations of Minkowski question mark function turned out to be important for studying fixed points of $?(x)$, that is the solutions of the equation 
\begin{equation}\label{fixed}
?(x)=x.
\end{equation} 
It is clear that $?(0)=0$, $?(\frac{1}{2})=\frac{1}{2}$, $?(1)=1$ and there exists at least one solution of \eqref{fixed} in the interval $(0,\frac{1}{2})$ and at least one in the interval $(\frac{1}{2},1)$. A famous conjecture claims that the Minkowski question mark function $?(x)$ has exactly five fixed points $0<x_1<\frac{1}{2}<x_2<1$. Moreover, there is only one irrational fixed point of $?(x)$ in the interval $(0,\frac{1}{2})$.\\
It is possible to approach this problem by studying iterations of $?(x)$, because of the following
\begin{prop}
If for every rational number $p/q\in [0,1]$ one has 
\begin{equation}\label{iterPROP}
\lim_{n\to\infty} f_n\left(\frac{p}{q}\right) = A, \text{  where }A\in \{0,\frac{1}{2},1\},
\end{equation}
 then there exist exactly five fixed points of $?(x)$ in the interval $[0,1]$.
\end{prop}
In fact, the opposite is also true: if the conjecture about fixed points is true, then \eqref{iterPROP} holds for every rational $p/q$. \\
Certain results about fixed points of the function $?(x)$ were recently obtained in \cite{Shulga}.
\vskip +1em
We will study iterations of $?(x)$ from a different perspective. Notice that for almost all $x_0\in[0,1]$ with respect to Lebesgue measure, one has
\begin{equation}\label{deriv} f_n'(x_0)=?'(x_0)\cdot ?'(x)|_{x=?(x_0)}\cdot ?'(x)|_{x=?(?(x_0))}\cdot\ldots\cdot?'(x)|_{x=\underbrace{?(?(...?}_\text{n-1 times}\bigl(x_0\bigl)\bigl)\bigl)}.\end{equation}

As $?(x)$ is a continuous strictly increasing map $[0,1]\to[0,1]$, so is $f_n(x)$. By the Lebesgue theorem, for almost all $x\in [0,1]$ function $f_n(x)$ has finite derivative.
It is well-known (see \cite{Salem}) that when derivative $?'(x)$ exists, it equals either $0$ of $+\infty$. Hence by the Lebesgue theorem for almost all $x$ one has $?'(x)= 0$, so the Minkowski question mark is a singular function.
However, for $n\ge2$ it is not known whether $f_n(x)$ is a singular function with a property analogous to $?(x)$.
For example, it is an open question if the derivative of $?(?(x))$ is either equal to $0$, $+\infty$ or does not exist.

Apart from obvious cases, it seems to be non-trivial to find explicit examples of $x_0$ such that $f_n'(x_0)=0$.
By obvious examples we mean all rational numbers (since \eqref{deriv} because of $?'(z)=0 $ and $?(z)\in\mathbb{Q}$ for any $ z\in \mathbb{Q}$), as well as some quadratic irrationals $x_0 = [a_1,\ldots,a_t,\ldots]$ satisfying 
\begin{equation}\label{kappa} 
\liminf_{t\to\infty}{\frac{a_1 +...+a_t}{t}} > \kappa_2=4.401^+, 
\end{equation}
where for $\kappa_2$  the exact formula
$$
\kappa_2 = \frac{ 4 L_5 - 5 L_4} { L_5 - L_4 }, \,\,\,\,\, L_j = \log \frac{j+\sqrt{j^2+4}}{2} - j \cdot \frac{\log2}{2}, \,\,\,\,\,\,\, j=4,5
$$
was discovered in \cite{MD}.
The explanation is that in \cite{MD} the following theorem was proved.\\
\begin{thmx}\label{moshd}
Let for irrational number $x_0 = [a_1,\ldots,a_n,\ldots]$ one has \eqref{kappa}.
Then $?'(x)$ exists and $?'(x)=0$.
\end{thmx}
 From Theorem A we see that if $x_0$ is a quadratic irrationality and formula \eqref{kappa} is satisfied then simultaneously $?(x_0)\in\mathbb{Q}$ and $?'(x_0)=0$ and so by \eqref{deriv} we have $f'_n(x_0)=0$.\\
	
In the present paper we introduce a set $M$ of irrational numbers, such that for all $x_0\in M$ and for all $n\in\mathbb{Z}_+$ one has $f_n'(x_0)=0$. 

\section{Main result}

Let $ A=(a_1,a_2,\ldots,a_k)$  be a sequence of positive integers of arbitrary length $k\ge 0$. When $k=0$ we suppose that $A$ is an empty sequence. For a non-empty sequence $A$ we consider the corresponding continued fraction
\begin{equation}\label{contfrac}[ A]=[a_1,a_2,\ldots,a_k] = \frac{1}{a_1+\frac{1}{a_2+\ldots+\frac{1}{a_k}}}.\end{equation}
For the sequence $ A=(a_1,a_2,\ldots,a_{k})$ we denote by \\
$$S_{ A} = a_1+a_2+\ldots+a_{k} $$
the sum of all its elements.\\
Given a sequence $ A=(a_1,a_2,\ldots,a_{k})$ we write
$d(A)=d(a_1,a_2,\ldots,a_k)= k$ for its length. \\
We denote by
$$
(A, b) = (a_1,\ldots, a_k, b)
$$
a sequence, obtained by concatenating sequence $A$ with an element $b\in\mathbb{Z}_+$.\\
Now let us consider a set of sequences $  A_i=(a_1^i,a_2^i,\ldots,a_{d(A_i)}^i)$.\\
By
$$
(A_1,A_2)=(a_1^1,a_2^1,\ldots,a_{d(A_1)}^1,a_1^2,a_2^2,\ldots,a_{d(A_2)}^2)
$$
we mean the sequence, obtained by concatenating sequences $A_1$ and $A_2$.\\
For the continued fraction of the form $x_0=[ A_1,\tau_1, A_2,\tau_2,\ldots]$, where $\tau_i \in\mathbb{Z}_+$, we use the following notation for the sum of all partial quotients up to ${A}_k$. 
$$  \sigma_{{A_k}} =  S_{A_1} + \tau_1 + S_{A_2}+\ldots+\tau_{k-1}+S_{A_k}.$$
\\
In the present paper we consider a special set $M$ of irrational numbers. We take arbitrary sequences
$A_i=(a_1^i,a_2^i,\ldots,a_{d(A_i)}^i)$ and consider continued fraction
\begin{equation}\label{formcont}
x=[A_1,\tau_1,A_2,\tau_2,\ldots]
\end{equation}
such that
\begin{equation}\label{tau_k}
\text{for all }k\in\mathbb{Z}_+\,\,\,\,\, \tau_k =  \sigma_{{A_k}}  + s_k, \,\,\, \text{where  } s_k > (\kappa_2 -1 ) S_{ A_{k+1}} + \sigma_{ A_{k}} +2,\,\, s_k\in\mathbb{Z_+} .
\end{equation}
Our set $M$ consists of all $x$ of the form \eqref{formcont} constructed by the procedure described, that is 
$$M=\Bigl \{ x \text{    of the form \eqref{formcont} with     } \tau_k \text{     satisfying \eqref{tau_k}} \Bigl \} .$$
 Our main result is the following \\
\begin{theor}
For any $x_0 \in M$ and for any $n\in\mathbb{Z}_+$ we have $f_n'(x_0) = 0$.
\end{theor}
Our special set $M$ contains some well-known constants, for example, Cahen's constant $$C=\sum_{i \geqslant 0} \frac{ (-1)^i}{S_i -1}, $$
 where $(S_n)$ is \textit{Sylvester's sequence} 
 $$S_0 = 2, \text{ and }S_{n+1} = S_n^2 - S_n+1 \text{ for }n\geqslant 0$$
 Indeed, as it was shown in \cite{davidson}, Cahen's constant can be expressed as a continued fraction
 $$C = [1, q_0^2, q_1^2, q_2^2, \ldots ],
 \text{\,\,where\,\,\,} q_0=1,q_1=1, \text{\,and\,\,\,} q_{n+2}=q_n^2 q_{n+1} +q_n \text{\,\,for\,\,} n \geqslant 0.$$
 It is easy to prove that $q_n > 2^{2^{n-3}}$ for $n \geqslant 3$. Hence, if we set 
 $$A_1 = (1, q_0^2, q_1^2, q_2^2, q_3^2), \tau_i = q_{3+i}^2 \text{\,\,for\,\,} i \geqslant 1, \text{\,\,and\,\,} A_{i} \text{\,\,to be empty for\,\,} i \geqslant 2,$$
we will get that $\tau_k > 2\sigma_{A_k} +2,$ from which it follows that $C$ belongs to the set $M$.




\section{Notation and preliminaries}
Let $A=(a_1,\ldots,a_k)$. If $d(A)=k\ge1$, then by $\overleftarrow A,A^{-}$ and $A_{-}$ we denote sequences $(a_{k},\ldots,a_1), (a_1,\ldots,a_{k-1})$ and $(a_2,\ldots,a_{k})$ respectively.\\
By $\langle A \rangle$ we denote the denominator of the continued fraction $[ A]$.\\
Using this notation, we recall the following classical formula (see \cite{Graham}):
$$\langle X \rangle \langle Y \rangle \le \langle X,Y\rangle=\langle X \rangle \langle Y \rangle+ \langle X^{-}\rangle \langle Y_{-} \rangle = \langle X \rangle \langle Y \rangle \Bigl(1+[\overleftarrow X][Y]\Bigl)\le2\langle X \rangle \langle Y \rangle.$$
Using this inequality twice, we get
\begin{equation}\label{continuants}a \langle X \rangle  \langle Y \rangle \le \langle X,a,Y \rangle \le (a+2)\langle X \rangle  \langle Y \rangle . \end{equation}
The last inequality will be crucial for our proof.\\


One of the key observations for our argument is following lemma, introduced by Gayfulin and Shulga in \cite{Shulga}, which is a variation of so-called \textit{folding lemma} (for details see \cite{shallit}).\\
\begin{lem}
   \textup{(\cite{Shulga}, Lemma 3.2)}
   Let $s$ be an arbitrary nonnegative integer and
$$
?([a_1,a_2,\ldots,a_{n-1}])=[b_1,b_2,\ldots,b_k],\,\,\,b_k\neq1.
$$
Consider the number 
$$
\theta =[a_1,a_2,\ldots,a_{n-1}, a_n],\,\,\,\text{where}\,\,\,a_n = \sum\limits_{i=1}^{n-1}a_i+s, \,\,\, s\ge0.$$
Then
\begin{enumerate}
\item{If $n\equiv k(mod$  $2)$, then $?(\theta)=[b_1,b_2,\ldots,b_{k-1},b_k-1,1,2^{s+1}-1,b_k,\ldots,b_1].$}
\item{If $n\equiv k+1(mod$  $2)$, then $?(\theta)=[b_1,b_2,\ldots,b_k,2^{s+1}-1,1,b_k-1,b_{k-1},\ldots,b_1].$}
\end{enumerate}
\end{lem}
Next lemma is an obvious corollary from the rule of comparison of values of continued fractions.\\
\begin{lem}
 Let $\textbf{X}<\textbf{Y}$ and
$$ \textbf{X}=[a_1,\ldots,a_n,x,...],$$ 
$$\textbf{Y}=[a_1,\ldots,a_n,y,...]. $$
Then every  $\textbf{Z}$, such that $\textbf{X}<\textbf{Z}<\textbf{Y}$, has a continued fraction of the form  
$$\textbf{Z}= [a_1,\ldots, a_n,z,\ldots],$$ where $min(x,y)\le z \le max(x,y) .$
\end{lem}
From Lemma 1 and Lemma 2 we deduce the following statement.\\
\begin{foll}
Let 
$$
?([a_1,a_2,\ldots,a_{n-1}])=[b_1,b_2,\ldots,b_k],\,\,\,b_k\neq1.
$$
 and $a_n = \sum\limits_{i=1}^{n-1}a_i+s$, where $s\in\mathbb{Z}_+$.\\
Consider the number  $\gamma=[a_1,\ldots,a_n,c_1,\ldots,c_s]$. Then there exists a sequence of positive integers $(b_{k+2},\ldots,b_p)$ for which 

\begin{enumerate}
\item{If $n\equiv k(mod$  $2)$, then $?(\gamma) = [b_1,\ldots,b_{k-1},b_k-1,1,z,b_{k+2},\ldots,b_{p}].$}
\item{If $n\equiv k+1(mod$  $2)$, then $?(\gamma) = [b_1,\ldots,b_k,z,b_{k+2},\ldots,b_{p}].$}
\end{enumerate}
where $$2^{s+1}-1\le z \le 2^{s+2}-1.$$
\end{foll}
\begin{solve} We give a proof only in the case where $n$ is even and $k$ is odd. Other cases are quite similar. By Lemma 1 we have
$$?([a_1,\ldots,a_n])=[b_1,\ldots,b_k,2^{s+1}-1,1,b_k-1,\ldots,b_1],$$ 
$$?([a_1,\ldots,a_n+1])=[b_1,\ldots,b_k,2^{s+2}-1,1,b_k-1,\ldots,b_1] .$$ 
Then for the number $\gamma=[a_1,\ldots,a_n,c_1,\ldots,c_s]$ we have 
$$[a_1,\ldots,a_n]<\gamma<[a_1,\ldots,a_n+1]$$\\
and so 
$$?([a_1,\ldots,a_n])<?(\gamma)<?([a_1,\ldots,a_n+1]).$$ \\
Then by Lemma 2 we get
$$?(\gamma) = [b_1,\ldots,b_k,z,b_{k+2},\ldots,b_{p}], \text{   where  } 2^{s+1}-1\le z \le 2^{s+2}-1.$$
\end{solve}

\section{Proof of the main result}
\begin{solve}
To prove Theorem 1 it is enough to verify two propositions:
$$\text{1) If }x_0\in\text{ M, then }?'(x_0)=0.$$
$$\text{2) If }x_0\in\text{ M, then }?(x_0)\in \text{ M}.$$

First, we will find the structure of the continued fraction for $?(x_0)$.\\
Let us define ${B_1}$ as a sequence of partial quotients of the image of Minkowski question mark function of $[ A_1]$ (which we can assume to be non-empty), that is\\
$$?([{A_1}]) = [{B_1}]. $$
If $ A_2$ is non-empty, we apply Corollary 1 setting
 $$(a_1,\ldots,a_{n-1}) = {A_1},$$
$$(b_1,\ldots,b_k) = {B_1}, $$
$$a_n = \tau_1=S_{{A_1}}+s_1, $$
$$(c_1,\ldots,c_s) = {A_2}.$$
Then for the number $\gamma= [ A_1,\tau_1, A_2]$ we have $?(\gamma)=[ B_1,z_1, B_2]$, where 
$$2^{s_1+1}-1\le z_{1} \le 2^{s_1+2}-1$$
 and by $ B_2$ we denoted a sequence $(b_{k+2},\ldots,b_p)$ from Corollary 1.\\
In order to define $ B_k$ with $k\ge3$ and $z_r$ with $r\ge2$ inductively, notice that by definition of the set $M$, for every $k\in\mathbb{Z}_+$, partial quotient $\tau_k$ is greater than the sum of all previous ones by $s_k$. When $ A_i$ is non-empty we can apply Corollary 1 with
$$(a_1,\ldots,a_{n-1}) = ( A_1,\tau_1, A_2,\ldots,\tau_{i-2}, A_{i-1}),$$ 
$$(b_1,\ldots,b_k)=( B_1,z_1, B_2,..., B_{i-1}),$$
$$a_n=\tau_{i-1},$$
 $$(c_1,\ldots,c_s) =  A_i .$$\\
We have $?([ A_1,\tau_1, A_2,\ldots,\tau_{i-2}, A_{i-1}])=[ B_1,z_1, B_2,...,z_{i-2}, B_{i-1}]$.  
Consider the number 
$$\gamma=[ A_1,\tau_1, A_2,\ldots,\tau_{i-2}, A_{i-1},\tau_{i-1}, A_i ].$$
Then by Corollary 1 we get
\begin{equation}\label{gamma123}
?(\gamma)=?([ A_1,\tau_1, A_2,\ldots,\tau_{i-2}, A_{i-1},\tau_{i-1}, A_i ])=[ B_1,z_1, B_2,..., B_{i-1},z_{i-1}, B_i],
\end{equation}
 where $2^{s_{i-1}+1}-1\le z_{i-1} \le 2^{s_{i-1}+2}-1$, and by $ B_i$ we denoted a sequence $(b_{k+2},\ldots,b_p)$ from Corollary 1.\\
If, for example, $ A_k$ is empty, then by Lemma 1

 $$
 ?([ A_1,\ldots,  A_{k-1}, \tau_{k-1}]) = [  B_1,\ldots,   B_{k-1}, z_{k-1},  B_k], \text{where, depending on parity,}
 $$
 \begin{equation}\label{empty}
    B_k = (1, b_{d(B_{k-1})}^{k-1}-1,  (\overleftarrow B_{k-1})_{-},\ldots, \overleftarrow B_1) \text{  or  }  B_k=( b_{d(B_{k-1})-1}^{k-1}+1,  (\overleftarrow B_{k-1})_{--},\ldots, \overleftarrow B_1) \end{equation}
  
  As we can see, even if $ A_k$ is an empty sequence, $ B_k$ is well-defined and is non-empty, so if the next sequence $ A_{k+1}$ is non-empty, we can continue the procedure for non-empty sequences (by Corollary 2). If it is empty, then we continue by Lemma 1.

 By continuing this procedure, we conclude that
 $$?(x_0)=?([ A_1,\tau_1, A_2,\ldots,\tau_{i-2}, A_{i-1},\tau_{i-1}, A_i,\ldots ]) = [ B_1,z_1, B_2,\ldots, B_{i-1},z_{i-1}, B_i,\ldots]$$
 with 
 \begin{equation}\label{z_i}
 2^{s_{i}+1}-1\le z_{i} \le 2^{s_{i}+2}-1
 \end{equation} for every $i\in\mathbb{Z}_+$, for every $x_0\in M$, that is of the form \eqref{formcont}.\\
 \vskip +1em
Now we will show that $?'(x_0)=0$. For the local convenience let us represent partial quotients of the continued fraction $[ A_1,\tau_1, A_2,\ldots]$ in a form 
$$x_0=[ A_1,\tau_1, A_2,\ldots]=[m_1,\ldots,m_t,\ldots].$$

By Theorem \ref{moshd}, to secure $?'(x_0)=0$ it is enough to show that 
\begin{equation}\label{derivM}
\liminf_{t\to\infty}{\frac{m_1 +\ldots+m_t}{t}} > \kappa_2.
\end{equation}
To estimate this limit from below, we make use of the following construction. Let us consider the continued fraction $$[p_1,\ldots,p_t,\ldots]=[\underbrace{1,\ldots,1}_{d(A_1)},\tau_1,\underbrace{1,\ldots,1}_{d(A_2)},\tau_2,\underbrace{1,\ldots,1}_{d(A_3)},\tau_3,\ldots],$$
where we basically took continued fraction of $x_0$ and replaced every partial quotient of $A_i$ with 1 for all $i$.
For any $t$ we have
\begin{equation}\label{ineqMP}
\frac{m_1 +\ldots+m_t}{t} \ge \frac{p_1 +\ldots+p_t}{t} .
\end{equation}
Now if we consider $\liminf$ of the right-hand side, it is clear that
\begin{equation}\label{numberP}
\liminf_{t\to\infty}{\frac{p_1 +\ldots+p_t}{t}} = \liminf_{j\to\infty}{\frac{\sum\limits_{i=1}^j d( A_i)+\sum\limits_{i=1}^{j-1} \tau_i}{   \sum\limits_{i=1}^j d( A_i) + (j-1)    }}  ,
\end{equation}
because we can calculate $\liminf$ of the sequence $\frac{p_1 +\ldots+p_t}{t},\, t=1,2,3,\ldots$ on the left-hand side of \eqref{numberP} by considering the subsequence 
$$t_j = d(A_1)+\ldots+d(A_j)+(j-1),\,\,\,\,\,\,\, j=1,2,3,\ldots$$ which corresponds to the expression on the right-hand side of \eqref{numberP}.\\
Now recall that $x_0\in M$, so we have \eqref{tau_k} and we can bound the right-hand side of \eqref{numberP} as 

$$\frac{\sum\limits_{i=1}^j d( A_i)+\sum\limits_{i=1}^{j-1} \tau_i}{   \sum\limits_{i=1}^j d( A_i) + (j-1)    } >  
\frac{\sum\limits_{i=1}^j d( A_i)+\sum\limits_{i=1}^{j-1} (2\sigma_{A_i}+(\kappa_2-1)S_{A_{i+1}}+2) }{   \sum\limits_{i=1}^j d( A_i) + (j-1)    } > \frac{\kappa_2 \sum\limits_{i=2}^j d( A_i)+\sum\limits_{i=1}^{j-1} (2\sigma_{A_i}+2) + d(A_1) }{   \sum\limits_{i=1}^j d( A_i) + (j-1)    }   $$
$$ > \frac{\kappa_2 \left( \sum\limits_{i=1}^j d( A_i) + (j-1) \right) +\sum\limits_{i=1}^{j-1} (2\sigma_{A_i}+2) -(\kappa_2-1)d(A_1)-\kappa_2(j-1)}{   \sum\limits_{i=1}^j d( A_i) + (j-1)    } $$
$$> \kappa_2 + \frac{\sum\limits_{i=1}^{j-1} (2\sigma_{A_i}+2) -(\kappa_2-1)d(A_1)-\kappa_2(j-1)}{   \sum\limits_{i=1}^j d( A_i) + (j-1)    } . $$
Very last fraction is greater than $0$ for $j$ large enough, because every $\sigma_i$ contains $S_{A_1} \ge d(A_1)$ as a term and $\sigma_{A_i}\ge 2^i$ for every $i\in\mathbb{Z}_+$.



Hence by \eqref{numberP}, \eqref{ineqMP} and \eqref{derivM} we get the desired equality $?'(x_0)=0.$

\vskip +1em
Now our goal is to prove that $?(x_0)\in M$.\\
We have to make sure that from  \eqref{tau_k} if follows that
\begin{equation}\label{g_k}
g_k := z_k - \sigma_{ B_k} > (\kappa_2 -1 ) S_{ B_{k+1}} + \sigma_{ B_{k}} +2
\end{equation}
for every $k\in\mathbb{Z}_+$.\\
From \eqref{g_k} one can see that $?(x_0)\in M$ and this is just what we need. In the rest of Section 4 we verify \eqref{g_k}. To verify it, we need to bound $z_k$, $\sigma_{ B_k}$ and $S_{ B_{k+1}}$.

First, we know that 
\begin{equation}\label{z_k}
z_k > 2^{s_k}.
\end{equation} 
To estimate $\sigma_{ B_k}$ and $S_{ B_{k+1}}$ we notice that if $\frac{p_j}{q_j}=[a_1,\ldots,a_j]$, then
$$
j \le a_1+\ldots+a_j \le q_j = \langle a_1,\ldots, a_j \rangle.
$$
This means that 
$$\sigma_{ B_k}< \langle B_1,\tau_1,\ldots,\tau_{k-1},B_k \rangle $$ and 
\begin{equation}\label{s_b}
S_{ B_{k+1}}< \langle B_{k+1} \rangle. 
\end{equation}
By formula \eqref{rational} we have $$?([ A_1,\tau_1, A_2,\ldots,\tau_{k-1}, A_k])=[ B_1,z_1, B_2,\ldots,z_{k-1}, B_k].$$ 
So
\begin{equation}\label{sigma_b}
 \sigma_{ B_k} = S_{B_1}+z_1+S_{B_2}+\ldots+z_{k-1}+S_{B_k}\le \langle B_1,\tau_1,\ldots,\tau_{k-1},B_k \rangle = 2^{\sigma_{ A_k}-1 }. 
 \end{equation}
Finally, we need to bound $\langle B_k \rangle$ with 
$k\in\mathbb{Z}_+$:\\
From \eqref{rational} we know that $\langle B_1 \rangle = 2^{S_{ A_1}-1}$. Next, to estimate $\langle B_i \rangle$ with $i\ge2$ we apply \eqref{continuants} with
$$X=( B_1,z_1, B_2,..., B_{i-1}),\,\, a = z_{i-1},\,\, Y =  B_i.$$
So we get the inequality
$$  z_{i-1}  \langle  B_1,z_1, B_2,..., B_{i-1} \rangle \langle  B_i \rangle \le \langle  B_1,z_1, B_2,..., B_{i-1},z_{i-1}, B_i \rangle \le (z_{i-1}+2)  \langle  B_1,z_1, B_2,..., B_{i-1} \rangle \langle  B_i \rangle $$
from which we deduce
\begin{equation}\label{contB}  \frac{\langle  B_1,z_1, B_2,..., B_{i-1},z_{i-1}, B_i \rangle}{(z_{i-1}+2)  \langle  B_1,z_1, B_2,..., B_{i-1} \rangle} \le  \langle  B_i \rangle \le  \frac{\langle  B_1,z_1, B_2,..., B_{i-1},z_{i-1}, B_i \rangle}{z_{i-1}  \langle  B_1,z_1, B_2,..., B_{i-1} \rangle}. \end{equation}
From \eqref{rational} we see that
$$\langle  B_1,z_1,  B_2,..., B_{i-1} \rangle = 2^{\sigma_{ A_{i-1}}-1},\,\,\,\langle  B_1,z_1, B_2,..., B_{i-1},z_{i-1}, B_i \rangle = 2^{\sigma_{ A_{i}}-1}.$$
By \eqref{z_i} we have
$$2^{s_{i-1}+1}-1\le z_{i-1} \le 2^{s_{i-1}+2}-1.$$
Recall that by \eqref{tau_k} we have
$$\tau_{i-1}=\sigma_{ A_{i-1}} + s_{i-1}.$$
Last four formulas enable us to write \eqref{contB} as 
\begin{equation}\label{b_i}
 2^{\sigma_{ A_{i-1}}+S_{ A_i}-3} \le\frac{2^{\tau_{i-1}+ A_i}}{2^{s_{i-1}+2}+1}\le  \langle  B_i \rangle \le \frac{2^{\tau_{i-1}+ A_i}}{2^{s_{i-1}+1}-1} \le 2^{\sigma_{ A_{i-1}}+S_{ A_i}}.
 \end{equation}
Note that both upper and lower bounds do not depend on $s_{i-1}$.\\
From \eqref{s_b} and upper bound from \eqref{b_i} we have
\begin{equation}\label{goga}
S_{B_{k+1}} \le \langle B_{k+1} \rangle  \le 2^{\sigma_{ A_{k}}+S_{ A_{k+1}}}. 
\end{equation}
Now we can apply bound \eqref{sigma_b} on $\sigma_{B_k}$, bound \eqref{goga} on $S_{B_{k+1}}$ as well as bound \eqref{z_k} on $z_k$ together with \eqref{tau_k} to get inequality \eqref{g_k}. 

Indeed, by \eqref{z_k} and \eqref{sigma_b} we have
\begin{equation}\label{last1}g_k \ge 2^{s_k} - 2^{\sigma_{A_k}-1}>2^{(\kappa_2 -1 ) S_{ A_{k+1}} + \sigma_{ A_{k}}+2} -2^{\sigma_{A_k}-1} .
\end{equation}
Meanwhile, we can use \eqref{s_b} and \eqref{sigma_b} to estimate the right-hand side of \eqref{g_k} from above as follows
\begin{equation}\label{last2}   (\kappa_2 -1 ) S_{ B_{k+1}} + \sigma_{ B_{k}} +2 < (\kappa_2 -1 )2^{\sigma_{ A_{k}}+S_{ A_{k+1}}} + 2^{\sigma_{A_k}-1}+2. 
\end{equation}

When $ S_{ A_{k+1}} \ge 1$ (i.e. for any non-empty sequence $ A_{k+1}$ ), we can see that right-hand side of \eqref{last1} is greater than right-hand side of \eqref{last2}.\\
Indeed, consider the inequality

$$2^{(\kappa_2 -1 ) S_{ A_{k+1}} + \sigma_{ A_{k}}} >(\kappa_2-1) 2^{\sigma_{ A_{k}}+S_{ A_{k+1}} } + 2^{\sigma_{ A_k}} +2$$
or, dividing both sides by $2^{\sigma_{ A_k}}$,
$$
2^{(\kappa_2 -1 ) S_{ A_{k+1}} } >(\kappa_2-1) 2^{S_{ A_{k+1}} } + 1 + \frac{2}{2^{\sigma_{ A_k}}},
$$
which is true under the condition $ S_{ A_{k+1}} \ge 1$. \\
To summarize, we verified \eqref{g_k} when $A_{k+1}$ is non-empty.
\vskip +1em
If $ A_{k+1}$ is an empty sequence, then from \eqref{empty} we have
$$ S_{ B_{k+1} }= \sigma_{B_{k}}. $$
From this we get that the right-hand side of \eqref{g_k} can be estimated from above as
$$(\kappa_2 -1 ) S_{ B_{k+1}} + \sigma_{ B_{k}} +2= \kappa_2 \sigma_{ B_{k}}+2\le \kappa_2 2^{\sigma_{ A_k}-1 }+2$$
(here in the last inequality we take into account \eqref{sigma_b}).\\
In our case when $A_{k+1}$ is an empty sequence, we can improve bound \eqref{last1} by noticing that $S_{A_{k+1}}=0$, so
$$g_k >2^{\sigma_{ A_{k}}+2} -2^{\sigma_{A_k}-1} .$$
Now we can see that 
$$g_k >2^{\sigma_{ A_{k}}+2} -2^{\sigma_{A_k}-1} >\kappa_2 2^{\sigma_{ A_k}-1 }+2 \ge (\kappa_2 -1 ) S_{ B_{k+1}} + \sigma_{ B_{k}} +2$$
and this is exactly what we needed. So in the case $A_{k+1}$ is empty, we have verified \eqref{g_k} also.\\

We have shown that \eqref{g_k} is true for all $k\in\mathbb{Z}_+$,
hence $?(x_0)\in M$.
\end{solve}
\vskip +1em
\textbf{Remark.} In the definition of the set $M$ we could have assume a weaker condition than 
$$
s_k > (\kappa_2 -1 ) S_{ A_{k+1}} + \sigma_{ A_{k}} +2.
$$
 For example, by estimating $\sigma_{ B_{k}}$ more carefully. We can do better if we consider this sum as
$$ 
\sigma_{ B_k} = \sum\limits_{i=1}^{k-1} z_i +\sum\limits_{i=1}^{k} S_{ B_i}
$$
and then bound every term separately. For the sake of readability we did not do it.

\section*{Acknowledgements}

The author thanks Nikolay Moshchevitin for careful reading and helpful comments.

N. Shulga is a scholarship holder of "BASIS" Foundation for Development of Theoretical Physics and Mathematics and is supported in part by the Moebius Contest Foundation for Young Scientists.

\end{document}